\documentclass[12pt]{article}
\begin{document}
\begin{center}
\large{\bf{A von Neumann theorem for uniformly distributed sequences of
partitions}}
\bigskip

{\large{\bf{Ingrid Carbone and Aljo\v{s}a Vol\v{c}i\v{c}}}}

Universit\`a della Calabria
\end{center}
\bigskip

\begin{abstract}
In this paper we consider permutations of sequences of partitions, obtaining a result which parallels von Neumann's theorem on permutations of dense sequences and uniformly distributed sequences of points.
\end{abstract}

\section {\bf Introduction} 
\bigskip 

S. Kakutani studied in [K] the notion of {\it uniformly distributed sequences of partitions} of the interval $[0,1]$. He introduced the following construction.
Fix a number $\alpha \in ]0,1[$. If $\pi$ is any partition of  $[0,1]$, its $\alpha$-refinement, denoted by $\alpha \pi$, is obtained subdividing the longest interval(s) of $\pi$ in proportion $\alpha/(1-\alpha)$. By $\alpha^n \pi$ we denote the $\alpha$-refinement of $\alpha^{n-1}\pi$.

Let $\omega=\{[0,1]\}$ be the trivial partition of $[0,1]$. The sequence $\{\alpha^n \omega\}$ will be called the Kakutani $\alpha$-sequence.
\bigskip

\noindent {\bf Definition 1.1} Given a sequence of partitions $\{\pi_n\}$ of $[0,1]$, with 
$$\pi_n = \{ [t_{i-1}^n , t_i^n], 1 \le i \le k(n)\}\,,$$ 
we say that it is {\it uniformly distributed (u.d.)}  if for any continuous function $f$ on $[0,1]$ we have 
$$\lim_{n\rightarrow \infty}
\frac{1}{k(n)}\sum_{i=1}^{k(n)} f(t_i^n)= \int_0^1 f(t)\,d t\,\eqno(1)$$ 
or, equivalently, if the sequence of discrete measures 
$$\frac{1}{k(n)}\sum_{i=1}^{k(n)} \delta_{t_i^n} \eqno (2)$$
converges weakly to the Lebesgue measure $\lambda$ restricted to $[0,1]$. Here $
\delta_{t}$ is the Dirac measure concentrated in $t$.
\smallskip

We can now state Kakutani's result (see [K] for more details).
\bigskip

\noindent {\bf Theorem 1.2} {\it The sequence $\{\alpha^n \omega\}$ is uniformly distributed.}
\bigskip

This was a partial answer to the following question posed by the physicist H. Araki, which regarded {\it random} splittings of the interval $[0,1]$. Let $X_1$ be choosen randomly (with respect to uniform distribution) on $[0,1]$ and then, inductively, once $X_1$, $X_2$, $\dots$, $X_n$ have been chosen, let $X_{n+1}$ be a point chosen randomly (and uniformly) in the largest of the $n+1$ intervals determined by the previous $n$ points. Can we conclude that the associated sequence of cumulative distribution functions converges uniformly, with probability $1$, to the distribution function of the uniform random variable on $[0,1]$? By a theorem due to Polya, this is equivalent to the condition expressed by (1). 

This question has been studied in [vZ] and [L] and later on in [PvZ].

Note that in the probabilistic setting we may neglect the possibility that the partition obtained in the $n$-th step has more than one interval of maximal length, since this event has probability zero. This is not the case for Kakutani's construction, since for every $\alpha$ the partition $\{\alpha^n \omega\}$ has, for infinitely many values of $n$, more than one interval of maximal length.

In [ChV] the notion of (deterministic)  u.d. sequences of partitions has been extended to probability spaces on complete metric spaces. In [CV] Kaktuani's splitting procedure has been extended to higher dimension. 

Recently some new results revived the interest for the subject. In [V] we introduced the concept of $\rho$-refinement of a partition $\pi$, which generalizes Kakutani $\alpha$ - sequence, and we proved that the sequence $\{\rho^n \omega\}$ is also u.d.  

We also investigated in [V] the connections of the theory of u.d. sequences of partitions to the well-established theory of uniformly distributed (u.d.) sequences of points, showing how is it possible to associate (many) u.d. sequences of points to any u.d. sequence of partitions.

\smallskip
A sequence of points $\{x_n\}$ is u.d. on $[0,1]$ if for every continuous function $f$ on $[0,1]$ we have
$$\lim_{n\rightarrow \infty} \frac{1}{n}\sum_{i=1}^{n} f(x_i)= \int_0^1 f(t)\,d t\,.$$ 
\smallskip

Since our methods allow to construct new u.d. sequences of points, we think these connections are interesting in view of possible applications to the quasi-Monte Carlo method, which is in the last decades the main motivation for the study of u.d. sequences of points (see [N]).

We continue to compare the two theories in the present paper, where we are concerned with the property analogous to the one studied by von Neumann in [vN], where it is proved the following theorem.
\bigskip

 \noindent {\bf Theorem 1.3} {\it If $\{x_n\}$ is a dense sequence of points in $[0,1]$, then there exists a rearrangement of these points, $\{x_{n_k}\}$, which is uniformly distributed.}
\bigskip

As we shall see, there is a remarkable difference between von Neumann's result and its extensions to Borel measures on $[0,1]$ and the corresponding result we obtain for sequences of partitions (Theorem 2.2). We shall comment on that in the last section of this paper.
\smallskip

One of the consequences of von Nemann's result is that there are many u.d. sequences of points.

Our purpose is analogous: we want to show that there are many u.d. sequences of partitions. Before we proceed, we need to define the concepts of permutation of a partition and that of density of a sequence of partitions.
\smallskip

\noindent {\bf Definition 1.4} Given a partition   $\pi = \{ [t_{i-1} , t_i], 1 \le i \le k\}$, we denote by $ l_i=t_i-t_{i-1}$ the lenght of its $i$-th interval. 
The {\it diameter} of $\pi$, denoted by diam$\,\pi$,  is equal to $\max _{1\le i\le k} l_i$.
\bigskip

\noindent {\bf Definition 1.5} Given a sequence of partitions $\{\pi_n\}$, we say that it is {\it dense} if 
$\lim_{n\rightarrow \infty}$ {\rm
diam}$\,\pi_n=0$.
\bigskip

\noindent {\bf Definition 1.6} If $\pi = \{ [t_{i-1} , t_i], 1 \le i \le k\}$ is a partition, its {\it permutation} is a partition $\pi' = \{ [s_{h-1} , s_h], 1 \le h \le k\}$ defined by the points  $s_h=\sum_{j=0}^h l_{i_j}\,,$ for $0\le h \le k$, where $l_0=0$ and $\{i_j\}$, for $1\le i \le k$, is a permutation of the indices $\{1, 2, \dots, k\}$.

 We will denote by $\pi !$ the set of all the permutations of $\pi$.
 
\bigskip
Of course, if $\pi$ splits $[0,1]$ in $k$ parts, $\pi$ has at most $k!$ permutations:  if $\pi$ has two or more intervals of the same length, the number of permutations is smaller than $k!$ . The extreme case is when all the intervals of $\pi$ have the same lenght. Then $\pi!$ coincides with the singleton \{$\pi\}$.

\smallskip
In the sequel we shall use the well known fact that uniform distribution of partitions can be tested by means of suitable families of functions and sets (see, for instance, [KN], Section 1 of Chapter 1 and for more general probability spaces, [DT], Section 2.1).
\bigskip

 \noindent {\bf Definition 1.7} We say that a family $\cal F$ of fuctions defined on $[0,1]$ is {\it determining} if a sequence of partitions $\{\pi_n\}$ is uniformly distributed if and only if (1) holds restricted to all $f \in \cal F$.
 \bigskip
 
It is very simple to show that the family of all characteristic functions of dyadic intervals
$I_h^s=[\frac{h-1}{2^s},\frac{h}{2^s}]$,  for $s\in{I\!\!N}$ and $1\le h\le 2^s$, is a determining family, i.e. a sequence of partitions $\{\pi_n\}$ is u.d. if and only if
$$\lim_{n\rightarrow
\infty} 
\frac{1}{k(n)}\sum_{i=1}^{k(n)} \chi_{I_h^s}(t_i^n)=\frac{1}{2^s} $$
for all $s\in{I\!\!N}$ and $1\le h\le 2^s$.
\smallskip

To simplify notation, we shall use in the sequel the following notation:
$$ \pi_n (I_h^s)=\frac{1}{k(n)}\sum_{i=1}^{k(n)} \chi_{I_h^s}(t_i^n),\eqno (3)$$ 
denoting the measure defined in (2) by the same symbol $\pi_n$ which is used for the corresponding partition.

\section{\bf The main result}

We begin this section by showing that density is a necessary condition for the uniform distribution of a sequence of partitions.
\bigskip

\noindent {\bf Proposition 2.1} {\it Any  uniformly
distributed sequence of partitions $\{\pi_n\}$ is dense.}
\bigskip

\noindent {\it Proof.} Suppose the contrary. Then there exists an
$\varepsilon >0$ such that, for infinitely many indices, diam$\,\pi_n \ge
\varepsilon$. Denote by
$n_k$ these indices and select for each $n_k$ an interval $I_{n_k}$ belonging to
$\pi_{n_k}$ having length at least $\varepsilon$. 

Let $m$ be a
positive integer such that $\frac{1}{m}$ is smaller than $\varepsilon$ and
consider the points
$x_i=\frac{i}{2m}$, for $0\le i \le 2m$. Then each  $I_{n_k}$ contains at
least two of them and therefore it contains at least one of the intervals
$J_i=[x_{i-1}, x_i]$, for $1\le i \le 2m$. It follows that at least one of
these intervals $J_i$ (denote it by $J$) is contained in infinitely
many
$I_{n_k}$'s. Let now $f$ be a continuous function having non
vanishing integral, whose support is contained in $J$. Then the sequence
$$ A_n=
\frac{1}{k(n)}\sum_{i=1}^{k(n)} f(t_i^n)$$ does not converge to
$\int_0^1f(t)dt$ when $n$ tends to infinity, since $A_{n_k}=0$ for all $k$. This contradicts the uniform distribution of $\{\pi_n\}$.
${\bf\star}$
\smallskip

The density is necessary for uniform distribution, but it is of course not sufficient. However, we have the
following result, which is the main result of this paper.

\bigskip

\noindent {\bf Theorem 2.2} {\it If $\{\pi_n\}$ is a dense sequence of
partitions,
then there exists a sequence of partitions $\{\sigma_n\}$, with
$\sigma_n\in{\pi_n!}$, which is uniformly distributed.}
\bigskip

\noindent {\it Proof.} Let $\{\pi_n\}$  be a dense sequence of partitions, and let  $k(n)$ denote the number of intervals of $\pi_n$. 

Since the binary intervals  $I_h^s=[\frac{h-1}{2^s},\frac{h}{2^s}]$ are a determining family, the conclusion will be acheaved constructing, for each $n$, a permutation  $\sigma_n$
 of $\pi_n$ such that
$$\lim_{n\rightarrow \infty}  \sigma_n (I_h^s)=\frac{1}{2^s} $$ 
\smallskip

\noindent for any $s\in{I\!\!N}$ and any $1\le h\le 2^s$ (see Definition 1.7 and formula (3)).

For any $s\in{I\!\!N}$, there exists $n_{s}\in{I\!\!}N$ such that  diam$\,\pi_ n \le \frac{1}{4^{s}}$
 for all $n\ge n_{s}$.   
 Of course, if $k(n)$ is the number of intervals of $\pi_n$, we have that $k(n)\ge4^{s}$  for all $n\ge n_{s}$.
 We may select a subsequence $\{n_{s}\}$ so that $n_{{s}+1}>n_{s}$ for all $s \in I\!\!N$. 
 \smallskip

When $n<n_1$, we just take $\sigma_n= \pi_n$. Suppose now $s>1$ and let us construct, for each  $n$ such that $n_{s} \le n<n_{{s}+1}$ , a permutation $\sigma_n$ of $\pi_n$ such that  $\sigma_n(I_h^t)$ is close to $\frac{1}{2^t}$ for all $t\le s$, in a sense which will be made precise later. 

 First, order the intervals of the partition $\pi_n$ with respect to their increasing length. 

Then we shift to the right the first interval of $\pi_n$ so that its right endpoint is $1$ and we move correspondingly the block of all the other intervals to the left. In this way the number of the intervals contained in $I_1^1$ remains unchanged or decreases by one unit.
We repeat this procedure until we obtain a permutation of $\pi_n$, denoted by $\pi_n^1$, such that
$$\frac{k(n)}{2} -1 \le k(n) \pi_n^1(I_h^1) \le  \frac{k(n)}{2} +1,\,\,\,\,\,\,\,\, 1 \le h \le 2. \eqno (4)$$

\smallskip

Now we consider $t=2$ and the corresponding dyadic intervals $I_h^2$, with $1 \le h \le 4$. Let us denote by $J_1^1=\left[{1 \over 2} - \delta _1^1(n), {1 \over 2} + \tilde \delta _1^1(n) \right[\,\,$ the interval of $\pi_n^1$ containing $1 \over 2$. Of course, $\tilde \delta_1^1 (n)>0$ and $  \delta_1^1(n)+\tilde \delta_1^1 (n) \le {1 \over 4^s } $.
\smallskip

We repeat in $I_1^1 \setminus J_1^1$ and in $I_2^1 \setminus J_1^1$ the procedure used above in order to get a convenient permutation of $\pi_n^1$ . More precisely, we keep $J_1^1$ fixed and  order the two collections of intervals of the partition $\pi_n^1$ contained in $I_1^1 \setminus J_1^1$ and $I_2^1 \setminus J_1^1$ respectively, with respect to their increasing lenght,  and in each of them we repeat the procedure described above, shifting the first interval to the right and the rest of them to the left and continue this procedure until the intervals are distributed, up to one unit,  proportionally to the lenghts of $I_1^2$ and $I_2^2 \setminus J_1^1$, respectively, i.e. proportionally to
$$ \alpha _1^2(n)= \frac{ {\frac {1}{4}}}{{\frac {1}{2}}- \delta_1^1 (n)}\,\,\,\, {\rm in} \,\,\,I_1^2 \,\,\,\, {\rm and}\,\,\,\,
\alpha _2^2(n)= \frac{ {\frac {1}{4}} - \delta_1^1(n)}{{\frac {1}{2}}- \delta_1^1 (n)}\,\,\,\, {\rm in} \,\,\,I_2^2\setminus J_1^1.$$

In the same way we reorder all the intervals of $\pi_n^1$ contained in $I_2^1 \setminus J_1^1$ until they are distributed proportionally (up to one unit) to 
$$\alpha _3^2(n)= \frac{ {\frac {1}{4}} - \tilde \delta_1^1(n)}{{\frac {1}{2}}- \tilde \delta_1^1 (n)}\,\,\,\, {\rm in} \,\,\,I_3^2\setminus J_1^1\,\,\,\, {\rm and}\,\,\,\,
\alpha _4^2(n)= \frac{ {\frac {1}{4}}}{{\frac {1}{2}}- \tilde \delta_1^1 (n)}\,\,\,\, {\rm in} \,\,\,I_4^2.$$

Since $\delta_1^1(n) + \tilde \delta_1^1(n)\le \frac{1}{4^{s}}$, it is clear that all  the coefficients $\alpha_h^2(n)$ may be estimated in terms of a function of $s$ in the following way: 

$$\frac{ {\frac {1}{4}} - {\frac {1}{4^{s}}} }{{\frac {1}{2}}} \le \alpha_h^2(n) \le \frac{ {\frac {1}{4}}}{{\frac {1}{2}} - {\frac {1}{4^{s}}}}, \,\,\,\, 1 \le h \le 2^2. \eqno (5)$$
\smallskip
\noindent Since the lower and upper estimates in (5) do not depend on $h$, and we are assuming that $n_s \le n < n_{s+1}$, in the sequel we will write $\alpha_s^2$ instead of $\alpha_h^2(n)$, to simplify notation.
\smallskip

We denote by $\pi_n^2$ the permutation of $\pi_n^1$ constructed above (which fixes $J_1^1$). Note that  $\pi_n^2$ does not move intervals of $\pi_n^1$ from $I_1^1$ to  $I_2^1$ and from  $I_2^1$ to  $I_1^1$; therefore  for $ h=1, 2 $ we can write 
$$(k(n)\pi_n^1(I_1^1)-1) \alpha_s^2 -1 \le k(n) \pi_n^2(I_h^2) \le  k(n)\pi_n^1(I_1^1)) \alpha_s^2 +1 $$ and
$$k(n)\pi_n^1(I_2^1) \alpha_s^2 -1 \le k(n) \pi_n^2(I_h^2) \le  k(n)\pi_n^1(I_2^1) \alpha_s^2 +1$$
\bigskip
\noindent for $h=3, 4$.

If we put $\alpha_s^1=\frac{1}{2}$,  taking (4) into account  we get
$$k(n)\alpha_s^1 \alpha_s^2 -2 \alpha_s^2-1 \le k(n) \pi_n^2(I_h^2) \le  k(n) \alpha_s^1 \alpha_s^2 + \alpha_s^2+1,\,\,\,\,1 \le h \le 2^2. $$

\smallskip

Let us make one more step in order to indicate how to procede for any $t \le s$.

Denote by $J_1^2$ and $J_2^2$ the two intervals of $\pi_n^2$ containing, respectively, the points $1 \over 4$ and $3 \over 4$ and set  
$$J_1^2=\left [{1 \over 4} -\delta _1^2(n), {1 \over 4} + \tilde \delta _1^2(n) \right [  \,\,\,{\rm and}\,\,\,J_2^2= \left[{3 \over 4} -\delta _2^2(n), {3 \over 4} + \tilde \delta _2^2(n) \right[$$

\noindent with  $\tilde \delta _1^2(n), \tilde \delta _2^2(n) >0$. We also note that
 $\lambda (J_1^2)$ and $\lambda (J_2^2)$ are both bounded from above by $  1 \over 4^{s}$.
 \smallskip
 
Let us procedeed as before, but now in the four intervals  $I_h^2 \setminus J_1^1 \cup J_1^2 \cup J_2^2$ for $1 \le h \le 2^3$.  By reordering all the intervals of $\pi_n^2$ contained in each of them we obtain a new partition $\pi_n^3 \in \pi_n^2!$ (which keeps fixed $J_1^1, J_1^2$ and $J_2^2$) whose intervals are distributed proportionally (up to one unit) to certain coefficients $\alpha_h^3(n)$ (which are determined as before), all of which satisfy the inequalities
 
 $$\frac{ {\frac {1}{2^3}} - {\frac {1}{4^{s}}} }{{\frac {1}{2^2}}} \le \alpha_h^3(n) \le \frac{ {\frac {1}{2^3}}}{{\frac {1}{2^2}} - {\frac {2}{4^{s}}}}, \,\,\,\, 1 \le h \le 2^3. $$

\smallskip 

Here we note again that, as in formula (5), the estimates do not depend on $h$, and $n$ is fixed, so in the sequel we shall write $\alpha_s^3$ instead on $\alpha_h^3(n)$. 
\smallskip

We note, as in the previous step, that $\pi_n^3$ does not move intervals of $\pi_n^2$ among the dyadic intervals $I_h^2$, with $1 \le h \le 2^2$.

 We easily get
 $$k(n)\alpha_s^1 \alpha_s^2\alpha_s^3 -2 \alpha_s ^2\alpha_s^3 -2\alpha_s^3-1 \le k(n)\pi_n^3(I_h^3)  \le k(n)\alpha_s^1 \alpha_s^2\alpha_s^3 + \alpha_s ^2\alpha_s^3 +\alpha_s^3+1$$ for all $ 1\le h \le 2^3. $
 
 \smallskip
 
 If is clear now that the same procedure can be repeated for all $t\le s$ and that we can construct partitions $\pi_n^t \in \pi_n^{t-1}!$ such that 
 $$k(n) \prod_{j=1}^t \alpha_s^j- 2 \,\,\sum_{i=2}^t \,\,\prod_{j=i}^t \alpha_s^j-1\le k(n) \pi_n^t(I_h^t)\le 
 k(n) \prod_{j=1}^t \alpha_s^j+ \sum_{i=2}^t \,\,\prod_{j=i}^t \alpha_s^j+1, \eqno(6)$$ 
\smallskip
\noindent where  $\alpha_s^t$, with $2\le t \le s$, are the proportionality coefficients which satisfy

 $$\alpha_s^1=  \frac {1}{2} \,\,\,\,\,\,\, {\rm and} \,\,\,\,\,\,\frac{ {\frac {1}{2^t}} - {\frac {1}{4^{s}}}}{{\frac {1}{2^{t-1}}}} \le \alpha^t_s \le \frac{ {\frac {1}{2^t}}}{{\frac {1}{2^{t-1}}} - {\frac {2}{4^{s}}}}, \,\,\,\,\,\, 2 \le t \le s. \eqno (7)$$
 
 \smallskip
Moreover, (6) implies the following estimates
 $$ k(n) \left | \,\, \pi_n^t(I_h^t)- \prod_{j=1}^t \alpha_s^j\,\,\right | \le 2\,\, \sum_{i=2}^t \,\,\prod_{j=i}^t \alpha_s^j+1, \,\,\,\,\, 1 \le h \le 2^t, \,\,1 \le t \le s.\eqno (8)$$ 

 Now we observe that  the partition $\pi_n^t$ does not move the intervals of $\pi_n^{t-1} $ contained in each dyadic interval $I_h^{t-1}$, with $1 \le h \le 2^{t-1}$, to another dyadic interval $I_k^{t-1}$, with $1 \le k \le 2^{t-1}$, i.e. $\pi_n^{t-1}(I_h^{t-1})=\pi_n^t(I_h^{t-1})$ for all $t\le  s$, $1\le h\le 2^{t-1}$.

Let us put $\sigma_n= \pi_n^s$ and observe that the condition $n_{s}\le n <n_{s+1}$  implies that  $n \rightarrow \infty$ if and only if $s \rightarrow \infty$.

If we divide the terms of (8) by $k(n)$ and substitute $\pi_n^t$ by $\sigma_n$, we obtain  
$$|\sigma_n(I_h^t)- L_s(t) | \le R_{n,s}(t), \,\,\,\, 1\le h \le2^t,\,\,\,\,1 \le t \le s,\eqno (9)$$

\noindent where
$$L_s(t)=  \prod_{j=i}^t \alpha_s^j\,,$$

\noindent and
$$ R_{n,s}(t)= \frac{1}{k(n)} \left (2\,\, \sum_{i=2}^t \,\,\prod_{j=i}^t \alpha_s^j+1 \right).$$

We will provide now estimates for  $L_s(t) - \ \frac{1}{2^t}$ and for $R_{n,s}(t) $. To this purpose we observe that from (7) we get  for all $1 \le t \le s$

$$\frac{1}{2^t} \prod_{j=2}^t \left(1- {\frac {2^j}{4^{s}}} \right) \le L_s(t) \le \frac{1}{2^t} \frac{ 1}{\prod_{j=2}^t  \left (1 - {\frac {2^j}{4^{s}}} \right)}\,.$$
\smallskip

\noindent If we substitute all the terms in the brackets of the previous formula by the smallest of them, we get

$$\frac{1}{2^t}< \frac{1}{2^t} \left(1- {\frac {2^t}{4^{s}}} \right)^{t-1} \le L_s(t) \le \frac{1}{2^t}   \frac{ 1}{\left(1 - {\frac {2^t}{4^{s}}}\right)^{t-1}}\,.$$

We note that 
$${\left(1 - {\frac {2^t}{4^{s}}}\right)^{1-t}}\le \left(\frac{4^s}{4^s-2^s}\right)^s\,\,\, \mbox {for all} \,\,\,1 \le t \le s $$
and  elementary calculation shows that the sequence of these upper bounds, which depend only on $s$, tends to $1$ when $s$ tends to infinity. Then, for all $1\le t\le s$ we have
$$  0< L_s(t) - \frac{1}{2^t} \le     \frac{1}{2^t} \left(\frac{4^s}{4^s-2^s}\right)^s -  \frac{1}{2^t} \eqno (10)$$

\smallskip 
\noindent and, consequently, 
$L_s(t)- \frac{1}{2^t} \rightarrow 0$ when $s$ tends to infinity.

\smallskip

 Moreover, for $R_{n,s}(t)$ we have the following estimate:
 $$0<R_{n,s}(t) \le \frac{2}{k(n)}  \,\, \sum_{i=2}^t \,\,\prod_{j=i}^t \alpha^j  +\frac{1}{k(n)} \le \frac{2s-2}{4^{s}} L_s(s) +\frac{1}{4^s}  \eqno (11)$$
 for all $1\le t\le s $, therefore $R_{n,s}(t) $ is bounded by a sequence which tends to zero when $s$ tends to infinity.
 
 Substituting (10) and (11) in formula (9), we get 
 
 $$ \left |\sigma_n(I_h^t)-\frac{1}{2^t} \right |\le R_{n,s}(t) \,+L_s(t) - \frac{1}{2^t} \le  \frac{2s-2}{4^{s}} L_s(s) +\frac{1}{4^s}+ \frac{1}{2^t} \left[ \left(\frac{4^s}{4^s-2^s}\right)^s -  1 \right]. $$
 
 \smallskip
 Putting all things together, we can conclude therefore that
    $$\sigma_n(I_h^t)  \rightarrow \frac{1}{2^t}  \,\,\, \rm{when}\,\,\, n \rightarrow \infty \,\,\, \mbox{for all} \,\,\,t \in {I\!\!N}. \eqno \star$$

\section {\bf Conclusions} 

We have constructed the u.d. sequence $\{\sigma_n\}$ from a dense sequence $\{\pi_n\}$ by a precise algorithm which could be, in concrete situations, implemented step by step. 

On the other hand it would be interesting to have besides our result also a probabilistic statement expressed by the following
\bigskip

\noindent {\bf Conjecture} {\it If $\{\pi_n\}$ is a dense sequence of partitions and if each $\sigma_n$ is  taken at random from $\pi_n!$ for any $n\in I\!\!N$, then $\{\sigma_n\}$  is uniformly distributed with probability $1$.}

\bigskip
The (possible) confirmation of the correctness of this conjecture would not diminish the value of the result presented in this paper. 
If it is allowed to compare the questions we are treating with much more important ones, we all agree that  the Borel theorem on normal numbers does not diminish the interest for finding concrete examples of such numbers.

\bigskip

It is interesting to note also an essential difference between the theorem due to von Neumann and our result. It is well known  that the conclusion of Theorem 1.3 holds if the measure $\lambda$ is substituted by any Borel probability on $[0,1]$ (see for instance Section 4 of Chapter 2 and, for a more general setting, Section 2 of Chapter 3, of [KN] and the bibliography cited there).
\smallskip

This is not the case with our result, as it can be easily seen by taking a sequence of partitions $\{\pi_n\}$, where each $\pi_n$ splits $[0,1]$ in intervals of equal length. Then for any $n\in I\!\!N$ the set of permutations $\pi_n !$ is a singleton and $\lambda$ is the only possible limit.
\smallskip

This is of course an extreme case and it is not difficult to see by examples that in many cases, if $\{\pi_n\}$ is a dense sequence of partitions, the corresponding set $\cal M$ of all the probabilities which are limits of sequences $\{\sigma_n\}$, with $\sigma_n\in \pi_n !$, contains more than just $\lambda$.
\smallskip

This observation rises questions about the size of the set $\cal M$  (could it contain in some cases all the Borel measures?) and its geometric and topological properties (is it convex, is it closed?).
\bigskip

The conjecture and the questions about  $\cal M$ will be addressed in subsequent papers.
\newpage

\noindent {\bf References}
 \bigskip\bigskip
 
 \noindent  [ChV] F. Chersi,  A. Vol\v{c}i\v{c},  $\lambda$-equidistributed sequences of partitions and a theorem of the de Bruijn-Post type,  {\it Annali Mat. Pura Appl. (4) {\bf 162 } } (1992) 23-32.
\smallskip

\noindent  [CV] I. Carbone, A. Vol\v{c}i\v{c}, Kakutani splitting procedure in
higher dimension, {\it  Rend. Ist. Matem. Univ. Trieste} {\bf 39}  (2007) 119-126.
\smallskip

\noindent [DT] M. Drmota,  R. F. Tichy,  {\it Sequences, discrepancies and applications},  Lecture Notes in Mathematics {\bf 1651}, Springer Verlag, Berlin, 1997. 
\smallskip

\noindent [K] S. Kakutani, A problem on equidistribution on the unit
interval $[0,1]$, {\it  Measure theory
(Proc. Conf., Oberwolfach, 1975)},  pp. 369--375. {\it Lecture Notes in Math.} {\bf 541}, Springer, Berlin, 1976.
\smallskip

\noindent [KN] L. Kuipers, H. Niderreiter, {\it Uniform distribution of sequences.  Pure and Applied Matematics}. Wiley-Interscience, New York-London-Sidney, 1974.
\smallskip

\noindent [L1] J. C. Lootgieter, Sur la r\'epartition des suites de Kakutani, {\it C. R. Acad. Sci. Paris S\'er. A-B}  {\bf \,285} (1977), no. 5, A403-A406.
\smallskip

\noindent [L2] J. C. Lootgieter, Sur la r\'epartition des suites de Kakutani, {\it C. R. Acad. Sci. Paris SŽr. A-B} {\bf 286} (1978), no. 10, A459-A461

\noindent [N] H. Niderreiter, {\it Random number generation and quasi-Monte Carlo Methods}, CBMS-NSF Regional Conference Series in Applied Math., 1992.
\smallskip

\noindent [PvZ]  R. Pyke, W. R. van Zwet, Weak convergence results for the Kakutani interval splitting procedure, {\it Ann. Probability} {\bf \,32} (2004), no. 1A, 380-423.
\smallskip

\noindent [vZ] W.R.  van Zwet, A proof of Kakutani's conjecture on random subdivision of longest intervals, {\it Ann. Probability} {\bf 6} (1978), no. 1, 133-137.
\smallskip

\noindent [V] A. Vol\v{c}i\v{c},  A generalization of Kakutani's splitting procedure. {\it Submitted}
\smallskip

\noindent [vN]  J. von Neumann, Gleichm\"assig dichte Zahlenfolgen, {\it Mat.
Fiz. Lapok} {\bf 32}  (1925) 32-40.
\smallskip

\noindent [W] H. Weyl, \"Uber ein Problem aus dem Gebiete der diophantischen Approximationen, {\it Nach. Ges. Wiss. G\"ottingen, Math.-phys. Kl.} (1914), 234-244.

\end{document}